\documentclass[a4paper,11pt]{article} 
\usepackage{authblk}
\usepackage{fullpage}
\usepackage{graphicx}
\usepackage{xcolor}
\usepackage{algorithm}
\usepackage{algpseudocode}
\usepackage{amssymb}
\usepackage{amsthm}
\usepackage{amsmath}
\usepackage{mathrsfs}
\usepackage[mathcal]{euscript}
\usepackage{bm}
\usepackage{mathabx}
\usepackage[square,numbers,comma,sort&compress]{natbib} 
\usepackage[colorlinks,urlcolor=cobalt,citecolor=cobalt,linkcolor=cobalt,pdftex, pdfauthor = Lukas Exl]{hyperref}
\usepackage{tabularx}
\usepackage{setspace} 
\usepackage{wasysym}

\definecolor{airforce}{rgb}{0.16,0.32,0.75}
\definecolor{cobalt}{rgb}{0.0,0.28,0.67}

\newcommand{\tenV}{\bm{\mathcal{V}}}

\newcommand{\tenA}{\bm{\mathcal{A}}}
\newcommand{\tenX}{\bm{\mathcal{X}}}
\newcommand{\tenY}{\bm{\mathcal{Y}}}
\newcommand{\tenC}{\bm{\mathcal{C}}}
\newcommand{\tenB}{\bm{\mathcal{B}}}

\newcommand{\Mr}{{\mathscr{M}_\mathbf{r}}}
\newcommand{\Id}{{\mathbf{I}}}
\newcommand{\Xn}{{\mathbf{X}_{(n)}}}
\newcommand{\Vn}{{\mathbf{V}_n}}
\newcommand{\Vnt}{{\mathbf{V}^T_n}}
\newcommand{\Bn}{{\mathbf{B}_{n}}}
\newcommand{\Bnt}{{\mathbf{B}^T_{n}}}
\newcommand{\Bnp}{{\mathbf{B}^\prime_{n}}}
\newcommand{\Cn}{{\mathbf{C}_{(n)}}}
\newcommand{\dCn}{{\dot{\mathbf{C}}_{(n)}}}
\newcommand{\Cnt}{{\mathbf{C}^T_{(n)}}}
\newcommand{\Un}{{\mathbf{U}_n}}
\newcommand{\Unh}{{\mathbf{U}^h_n}}
\newcommand{\Uell}{{\mathbf{U}_\ell}}
\newcommand{\Uellt}{{\mathbf{U}^T_\ell}}
\newcommand{\dUell}{{\dot{\mathbf{U}}_\ell}}
\newcommand{\Unt}{{\mathbf{U}^T_n}}
\newcommand{\dUn}{{\dot{\mathbf{U}}_n}}

\newcommand{\deltaUn}{{{\Delta\mathbf{U}}_n}}
\newcommand{\deltaUk}{{{\Delta\mathbf{U}}_k}}
\newcommand{\Uk}{{\mathbf{U}_k}}
\newcommand{\Ukt}{{\mathbf{U}^T_k}}
\newcommand{\dUk}{{\dot{\mathbf{U}}_k}}

\newcommand{\Vk}{{\mathbf{V}_k}}
\newcommand{\Vkinv}{{\mathbf{V}^{-1}_k}}
\newcommand{\Qk}{{\mathbf{Q}_k}}
\newcommand{\Qkt}{{\mathbf{Q}^T_k}}
\newcommand{\Rn}{{\mathbf{R}_n}}
\newcommand{\Rnt}{{\mathbf{R}^T_n}}
\newcommand{\tRn}{{\widetilde{\mathbf{R}}_n}}
\newcommand{\tRnt}{{\widetilde{\mathbf{R}}^T_n}}
\newcommand{\Tn}{{\mathbf{T}_n}}
\newcommand{\Tnt}{{\mathbf{T}^T_n}}
\newcommand{\Wn}{{\mathbf{W}_n}}
\newcommand{\Wnt}{{\mathbf{W}^T_n}}

\newcommand{\Pn}{{\mathbf{P}_n}}
\newcommand{\Pno}{{\mathbf{P}^{\bot}_n}}
\newcommand{\TYMr}{{\mathscr{T}_{\tenY} \mathscr{M}_\mathbf{r}}}
\newcommand{\subTYMr}{\widetilde{\mathscr{T}}_{\tenY} \mathscr{M}_\mathbf{r}}

\newcommand{\VIr}{{\mathcal{V}_{I,r}}}
\newcommand{\TUnVInrn}{{\mathscr{T}_{\Un} \mathcal{V}_{I_n,r_n}}}

\newcommand{\TUVIr}{{\mathscr{T}_{\mathbf{U}} \mathcal{V}_{I,r}}}
\newcommand{\RgUn}{{\mathscr{R}(\Un)}}
\newcommand{\RgUno}{{\mathscr{R}(\Un)}^\bot}
\newcommand{\RgU}{{\mathscr{R}(\mathbf{U})}}
\newcommand{\RgUo}{{\mathscr{R}(\mathbf{U})}^\bot}
\newcommand{\Ab}{{\mathbf{A}}}
\newcommand{\Bb}{{\mathbf{B}}}
\newcommand{\Rb}{{\mathbf{R}}}

\newcommand{\Tb}{{\mathbf{T}}}

\newtheorem{thm}{Theorem}

\title{\Large{\textbf{An optimization approach for dynamical Tucker tensor approximation}}}
\author[1,2]{Lukas Exl \thanks{\texttt{lukas.exl@univie.ac.at}}}
\affil[1]{\small Faculty of Mathematics, University of Vienna, 1090 Vienna, Austria.} 
\affil[2]{\small Institute for Analysis and Scientific Computing, Vienna UT, 1040 Vienna, Austria.} 

\begin{document}
\date{}
\maketitle

\noindent\textbf{Abstract.} An optimization-based approach for the Tucker tensor approximation of parameter-dependent data tensors and solutions of tensor differential equations with low Tucker rank is presented. 
The problem of updating the tensor decomposition is reformulated as fitting problem subject to the tangent space without relying on an orthogonality gauge condition. A discrete Euler scheme is 
established in an alternating least squares framework, where the quadratic subproblems reduce to trace optimization problems, that are shown to be 
explicitly solvable and accessible using SVD of small size. In the presence of small singular values, instability for larger ranks is reduced, since the method does not need the (pseudo) inverse 
of matricizations of the core tensor. Regularization of Tikhonov type can be used to compensate for the lack of uniqueness in the tangent space. The method is validated numerically  
and shown to be stable also for larger ranks in the case of small singular values of the core unfoldings. Higher order explicit integrators of Runge-Kutta type can be composed. \\



\section{Introduction}\label{intro}
Low-rank matrix and tensor approximation occurs in a variety of applications as model reduction approach \cite{grasedyck2013literature}. 
This paper concerns the problem of fitting parameter or time-dependent matrices and tensors with small (Tucker) rank within the Tucker tensor format \cite{kolda2009tensor}, that is, 
updating the factors and core tensor of the model instead of working with the tensor itself. 
{In a discrete-time setting for tensor streams, incremental fitting procedures for dimensionality reduction were introduced as dynamic higher order generalizations of principal component analysis \cite{sun2006beyond}. 
While this approach seeks to minimize a reconstruction error by a sequence of Tucker tensor approximations sharing same factors and is discrete in time, a different time continuous setting}
was recently discussed {by Koch and Lubich} \cite{koch2007dynamical,koch2010dynamical}, where differential equations for 
the factor matrices and core tensor of the Tucker approximation were derived. For a time-varying family of tensors $\tenA(t) \in \mathbb{R}^{I_1 \times I_2 \times \cdots \times I_d}$ and the approximation manifold $\Mr$ 
(Tucker tensors with rank $\mathbf{r}$) an alternative to the best-approximation problem 
\begin{align}\label{eqn:bestappr}
  \min_{\tenX(t) \in \Mr} \|\tenX(t) - \tenA(t) \|
\end{align}
is to considered a dynamical tensor approximation problem based on derivatives $\dot{\tenY}$ which lie in the tangent space $\TYMr$, i.e., 
\begin{align}\label{eqn:dynopt}
 \min_{\dot{\tenY}(t) \in \TYMr} \|\dot{\tenY}(t) - \dot{\tenA}(t) \|
\end{align}
with an initial condition $\tenY(0)$, typically a best-approximation in $\Mr$ of $\tenA(0)$. 
Note that for the choice $\dot{\tenA}(t) = F(t,\tenY(t))$ this becomes the defect approximation of a tensor differential equation $\dot{\tenY}(t) = F(t,\tenY(t))$. 
A prominent example of a tensor differential equation on $\Mr$ is the multiconfiguration time-dependent {Hartree} approach (MCTDH) 
of molecular quantum dynamics \cite{beck2000multiconfiguration,lubich2015time}, where a Schr\"odinger equation subject to $\Mr$ is solved. More precisely, the MCTDH ansatz takes the form 
\begin{align*}
  \mathrm{i}\,\dot{\psi} = \mathscr{H} \psi \quad {\text{subject to}} \quad \psi(\mathbf{x},t) = \sum_{j_1=1}^{r_1} \cdots \sum_{j_d=1}^{r_d} a_{j_1\cdots j_d}(t) \prod_{\kappa=1}^d \varphi_{j_\kappa}^{(\kappa)} (x_\kappa,t),
 \end{align*}
where the discretized constraint on the many-particle wave function $\psi$ leads to a Tucker tensor {manifold}. A time-dependent variational principle is used to extract differential equations for the coefficient tensor 
$\big(a_{j_1\cdots j_d}\big)$ and the factors $\varphi_{j_\kappa}^{(\kappa)}$. This procedure is equivalent to solving problem \eqref{eqn:dynopt} for $\dot{\tenA}(t) =  \mathscr{H} {\tenY(t)} $, with $\tenY$ denoting 
the discretized wave function in Tucker format, by applying the Galerkin condition 
\begin{align*}
 \langle {\mathrm{i}}\,\dot{\tenY} - \mathscr{H} \tenY, \tenV \rangle = 0 \quad \text{for all}\quad \tenV \in \subTYMr,
\end{align*}
for a gauged tangent space $\subTYMr$. In this sense, problem \eqref{eqn:dynopt} yields an initial value problem on $\Mr$ in form of a system of nonlinear differential equations for the factors and core tensor in the Tucker 
decomposition. This system is explicit if 
an \textit{orthogonality gauge condition} in $\TYMr$ is imposed. However, the arising equations involve the (pseudo) inverse of the core unfoldings, which might be close to singular. Consider for instance the two dimensional 
(low-rank matrix) case, where $\Mr$ consists of matrices  $\mathbf{Y} = \mathbf{U} \mathbf{C} \mathbf{V}^T \in \mathbb{R}^{I_1 \times I_2}$ with 
$\mathbf{U}\in \mathbb{R}^{I_1 \times r}$, $\mathbf{V}\in \mathbb{R}^{I_2 \times r}$ and $\mathbf{C} \in \mathbb{R}^{r \times r}$ nonsingular. The mentioned equations for the factors of the model take the form (omitting the argument $t$)
\cite{koch2007dynamical}
\begin{eqnarray}\label{eqn:dynlow2d}
\begin{aligned}
\dot{\mathbf{U}} = &\, (\Id - \mathbf{U}\mathbf{U}^T) \dot{\mathbf{A}} \mathbf{V} \mathbf{C}^{-1},\\
\dot{\mathbf{V}} = &\, (\Id - \mathbf{V}\mathbf{V}^T) \dot{\mathbf{A}} \mathbf{U} \mathbf{C}^{-T},\\
\dot{\mathbf{C}} = &\, \mathbf{U}^T \dot{\mathbf{A}} \mathbf{V}.
\end{aligned}
\end{eqnarray}
Nonzero singular values of $\mathbf{Y}$ are those of the core $\mathbf{C}$, which has to be inverted in the first two equations of \eqref{eqn:dynlow2d}. Hence, very small singular values lead to severe problems concerning stability of numerical explicit or 
implicit integrators \cite{kieri2016discretized}. Unfortunately, this is a very realistic case, since in practice the prescribed initial ranks should be chosen larger than the unknown effective ranks (overestimation) because 
underestimating necessarily leads to large {model errors}. Therefore, a numerical scheme for \eqref{eqn:dynopt} has to manage the rank deficient case. 
Recently, Lubich and Oseledets \cite{lubich2014projector} introduced an operator splitting approach for the low-rank matrix approximation, which behaves robust in the mentioned problematic case. In the meanwhile 
the method was extended to the tensor train format \cite{lubich2015time2,oseledets2011tensor} and applied to MCTDH \cite{lubich2015time} {with related analysis presented in \cite{lubich2017time}}. 
{Recently, also projected time-stepping methods for tensor differential equations on a low-rank tensor train manifold where established \cite{dolgov2012fast} and analyzed \cite{Kieri_V:2017}.}
{Moreover, we stress that in applications the equations \eqref{eqn:dynlow2d} (resp. higher dimensional versions) are often \textit{not} well-defined due to initial states, 
which exhibit rank(s) strictly smaller than $\mathbf{r}$, leading to singular values of core unfoldings that are exactly zero. As a consequence, the factor derivatives or core unfoldings have to be initially modified appropriately or, 
if feasible, the full tensor equations are evolved over a short time to generate a proper initial state \cite{beck2000multiconfiguration,nonnenmacher2008dynamical}. }\\ 
In this paper we break out in a different direction. We will develop an integration scheme using the optimization framework \eqref{eqn:dynopt} itself rather than a system of nonlinear differential equations. 
{As in previous approaches the approximation manifold is equipped with a full rank condition for the factors
while the rank deficient case of core unfoldings is tackled by regularizing a fitting problem which exhibits the resulting (local) error.}
For that purpose, it turns out to be beneficial to give up the gauge condition in the tangent space. Remarkably, the proposed scheme does not need any inversion {or} (pseudo) inverse of the core unfoldings (density matrices in MCTDH). 
The approach is motivated from static fitting problems of data tensors, where quasi best-approximation in the Tucker format can be derived very efficiently 
from \textit{alternating least squares} (ALS) or general nonlinear minimization of the norm of the residual \cite{kolda2009tensor,hackbusch2012tensor}. For instance, the key of the \textit{higher order orthogonal iteration} (HOOI) \cite{de2000best} is the 
reduction of the quadratic subproblems in the ALS scheme to the task of computing an orthonormal basis of a dominant eigenspace \cite{kokiopoulou2011trace}. 
{We stress that alternating linear schemes, where ALS is a special case, were successfully applied to optimization tasks subject to a low-rank tensor train manifold \cite{holtz2012alternating} and 
local convergence was proven in the case of convex functionals \cite{rohwedder2013local}.} 
Here, we establish a discrete Euler scheme {for the dynamical problem} by alternating least squares, where we impose orthonormal columns for the updated factor matrices as a constraint. As a consequence, the originally quadratic subproblems get linear with explicit solutions 
via SVD. The method does not require any pseudo inverse of unfoldings of the core. Hence, the approach is far less affected by small singular values, although the stability of the SVD solutions 
of the involved trace optimization problems can still be negatively affected in severe cases. It turns out that Tikhonov regularization reduces this problem while also being a compensation for the lack of uniqueness 
of the tangent space representation. The regularization only acts on the column space of the factor matrices and would be inactive if the conventional orthogonality gauge condition would have been imposed.\\
In the following section we introduce tensor notation and give some details to Tucker decompositions and the approximation manifold. In section \ref{sec:kochlubich} we review the approach of Koch and Lubich \cite{koch2010dynamical}. 
Section \ref{sec:alseuler} formulates the optimization framework and establishes the Euler scheme. Finally, we validate our approach by means of numerical examples.  
\section{Prerequisites}
In the following a brief description of tensors, Tucker decomposition and the manifold of rank-$(r_1,\hdots,r_d)$ tensors is given. The reader is referred to the review article \cite{kolda2009tensor} or the book \cite{hackbusch2012tensor} for 
further details and references therein.\\
We treat real-valued matrices and tensors, but emphasize that the complex case (as in quantum mechanics applications) does not pose any problems. 
Specifically, the derivation of the Euler scheme of Sec.~\ref{sec:alseuler} with its Theorem~\ref{thm:trace} generalizes to this case. 
\subsection{Notation}
We denote matrices by boldface capital letters, e.g., $\mathbf{A}$ and vectors by boldface lowercase letters, e.g., $\mathbf{a}$. A tensor is denoted by boldface Euler script letters, e.g., $\tenX$. 
Scalars are denoted by lowercase letters, e.g., $\alpha$.\\
We use the Frobenius inner product and norm for matrices, denoted by $\langle.,.\rangle$ and $\|.\|$, respectively. The \textit{norm of a tensor} 
$\tenX \in \mathbb{R}^{I_1 \times I_2 \times \cdots \times I_d}{, \, d \geq 2,}$ is the square root of 
the sum of the squares of all elements, or equivalently, the Euclidean norm of the \textit{vectorized} tensor
\begin{align}
 \|\tenX\|  = \sqrt{\mathrm{vec}(\tenX)^T \mathrm{vec}(\tenX)}.
\end{align}
The \textit{inner product of two tensors} $\tenX,\tenY \in \mathbb{R}^{I_1 \times I_2 \times \cdots \times I_d}$ is the sum of the product of their entries, i.e.,
\begin{align}
 \langle \tenX, \tenY \rangle = \mathrm{vec}(\tenX)^T \mathrm{vec}(\tenY).
\end{align}
There holds $\|\tenX\|^2 =  \langle \tenX, \tenX \rangle$.\\
Fibers of a tensor $\tenX \in \mathbb{R}^{I_1 \times I_2 \times \cdots \times I_d}$ are the generalizations of matrix rows and columns. A mode-$n$ fiber of $\tenX$ consists of the elements with all indices fixed except the $n$th index.  
The $n$th \textit{unfolding of a tensor} $\tenX \in \mathbb{R}^{I_1 \times I_2 \times \cdots \times I_d}$ is the matrix $\Xn \in \mathbb{R}^{I_n \times \prod_{j \neq n}I_j }$ that arranges the mode-$n$ fibers to be the 
columns of the resulting matrix. We denote the corresponding matricization as $[\tenX]_{(n)} = \Xn$. \\
The $n$-rank of a tensor $\tenX \in \mathbb{R}^{I_1 \times I_2 \times \cdots \times I_d}$ is 
\begin{align}
r_n = \mathrm{rank}(\Xn),
\end{align}
and the vector $\mathbf{r} = (r_1,\hdots,r_d)^T$ is called the \textit{Tucker rank} of $\tenX$.\\
The \textit{$n$-mode product} of a tensor $\tenX \in \mathbb{R}^{I_1 \times I_2 \times \cdots \times I_d}$ with a matrix $\mathbf{U} \in \mathbb{R}^{J_n \times I_n}$ 
results in a tensor
\begin{align}
 \tenX \times_n \mathbf{U} \in \mathbb{R}^{I_1 \times \cdots\times I_{n-1} \times J_n \times I_{n+1} \times \cdots  \times I_d},
\end{align}
which is defined as
\begin{align}
 [\tenX \times_n \mathbf{U} ]_{(n)} = \mathbf{U} \Xn.
\end{align}
There holds (with appropriate dimensions) 
\begin{align}
 (\tenX \times_n \mathbf{U}) \times_n \mathbf{V} = \tenX \times_n \mathbf{V}\mathbf{U}
\end{align}
and for $n \neq m$
\begin{align}
 \tenX \times_n \mathbf{U} \times_m \mathbf{W} = \tenX \times_m \mathbf{W} \times_n \mathbf{U}.
\end{align}
We will further denote the range of a matrix $\mathbf{U} \in \mathbb{R}^{J_n \times I_n}$ (column space) with $\RgU$ and the orthogonal complement with $\RgUo$. 

\subsection{Tucker decomposition and manifold of rank-$(r_1,\hdots,r_d)$ tensors}
A \textit{Tucker decomposition} \cite{Tucker1966} of a tensor $\tenX \in \mathbb{R}^{I_1 \times I_2 \times \cdots \times I_d}$ has the form
\begin{align}\label{eqn:tucker}
\tenX = \tenC \times_1 \mathbf{U}_1 \cdots \times_d \mathbf{U}_d =: \tenC \bigtimes_{k=1}^d \Uk,
\end{align}
where the \textit{core tensor} $\tenC \in \mathbb{R}^{r_1 \times r_2 \times \cdots \times r_d}$ and the \textit{factor matrices} $\Uk \in \mathbb{R}^{I_k \times r_k}$. The representation \eqref{eqn:tucker} 
is not unique: For $\Vk \in \mathbb{R}^{r_k \times r_k}$ nonsingular, modifying the core $\widetilde{\tenC} = \tenC \bigtimes_{k=1}^d \Vk$ 
together with modified factor matrices $\widetilde{\mathbf{U}}_k = \Uk \Vkinv$ gives the same tensor $\tenX = \widetilde{\tenC} \bigtimes_{k=1}^d \widetilde{\mathbf{U}}_k$.\\
In the following we further assume the factor matrices having orthonormal columns, that is 
\begin{align}
\Ukt \Uk = \Id.
\end{align}
Note that this does not lead to a unique Tucker decomposition, since $\widetilde{\mathbf{U}}_k = \Uk \Qk$ with $\Qk$ orthogonal matrices and $\widetilde{\tenC} = \tenC \bigtimes_{k=1}^d \Qkt$ gives the same tensor $\tenX = \widetilde{\tenC} \bigtimes_{k=1}^d \widetilde{\mathbf{U}}_k$.\\
A Tucker decomposition can be written as a linear combination of rank-1 tensors that are formed as the outer products of the columns of the factor matrices, that is
\begin{align}
\tenX = \sum_{j_1,\hdots,j_d} c_{j_1,\hdots,j_d} \, \textbf{u}_{j_1}^{(1)} \circ \cdots \circ \textbf{u}_{j_d}^{(d)}, 
\end{align}   
where $\tenC = (c_{j_1,\hdots,j_d})$. 
The $n$th unfolding of a Tucker tensor is given by
\begin{align}
\Xn = \Un \Cn (\mathbf{U}_d \otimes \cdots \otimes \mathbf{U}_{n+1} \otimes \mathbf{U}_{n-1} \otimes \cdots \otimes \mathbf{U}_1)^T =: \Un \Cn \bigotimes_{k\neq n} \Ukt,
\end{align}
where $\otimes$ denotes the \textit{Kronecker product}.\\

For given Tucker rank $\mathbf{r} = (r_1,\hdots,r_d)^T$ the set 
\begin{align}
\Mr = \{ \tenX \in \mathbb{R}^{I_1 \times I_2 \times \cdots \times I_d}: \, r_n = \mathrm{rank}(\Xn),\, n = 1,\hdots,d\}
\end{align}
is the manifold of  rank-$(r_1,\hdots,r_d)$ tensors. Every element $\tenX \in \Mr$ has a representation as a Tucker decomposition, where w.l.o.g we assume the factor matrices having orthonormal columns, that is $\Ukt \Uk = \Id$. 
We will use $\Mr$ as an approximation manifold with (typically) $r_n \ll I_n$. In the special case of $d=2$ elements of $\Mr$ are usually referred to as  \textit{low-rank matrices}, where $\mathbf{r} = (r,r)^T$. 
In this case, the conventional notation is $\mathbf{X} = \mathbf{U} \mathbf{C} \mathbf{V}^T \in \mathbb{R}^{I_1 \times I_2}$ with $\mathbf{U}\in \mathbb{R}^{I_1 \times r}$, $\mathbf{V}\in \mathbb{R}^{I_2 \times r}$ and $\mathbf{C} \in \mathbb{R}^{r \times r}$ nonsingular, but not necessarily diagonal as in the case of a (truncated) singular value decomposition (SVD).\\
For $\tenY \in \Mr$ the corresponding \textit{tangent space} is denoted with $\TYMr$, where $\dot{\tenY} \in \TYMr$ has the form 
\begin{align}\label{eqn:tanvec}
\dot{\tenY} = \dot{\tenC} \bigtimes_{n=1}^d \Un + \sum_{n=1}^d \tenC \times_n \dUn \bigtimes_{k \neq n} \Uk,
\end{align}    
with $\dot{\tenC} \in \mathbb{R}^{r_1 \times r_2 \times \cdots \times r_d}$ and $\dUn \in \TUnVInrn$, where 
\begin{align}
\TUVIr := \{\dot{\mathbf{U}} \in \mathbb{R}^{I \times r}:\, \dot{\mathbf{U}}^T \mathbf{U} + \mathbf{U}^T \dot{\mathbf{U}} = 0\}
\end{align}
is the tangent space of $\mathbf{U} \in \VIr$, the \textit{Stiefel manifold} of real $I \times r$ matrices with orthonormal columns. If the \textit{gauge condition} (orthogonality constraint)
\begin{align}\label{eqn:orthgauge}
 \Unt\dUn = 0 \quad n=1,\hdots,d
\end{align} 
is imposed, we denote the tangent space with $\subTYMr$. {We stress that the gauge leads to a unique parametrization of the tangent space, in fact, the $\RgU$-component of the variation can be absorbed by the 
first term in \eqref{eqn:tanvec} \cite{hackbusch2012tensor}. 
}
\section{Related approach of Koch and Lubich}\label{sec:kochlubich}
Dynamical approximation of time-dependent tensors $\tenA(t)$ were previously considered in \cite{koch2010dynamical} with help of an orthogonal projection on the tangent space $\subTYMr$. 
This can be accomplished by a Galerkin condition, also known as \textit{Dirac-Frenkel variational {principle}} \cite{beck2000multiconfiguration}, that is (omitting the argument $t$)
\begin{align}\label{eqn:DiracFrenkel}
 \langle \dot{\tenY}-\dot{\tenA}, \tenV \rangle = 0 \quad \text{for all}\quad \tenV \in \subTYMr.
\end{align}
The core and the factors of the decomposition of $\tenY$ fulfilling \eqref{eqn:DiracFrenkel} are explicitly given by a system of differential equations. We state here the associated theorem from \cite{koch2010dynamical}.
\begin{thm}[\cite{koch2010dynamical}]
 For a tensor $\tenY \in \Mr$ with $n$-mode factors having orthonormal columns the condition \eqref{eqn:DiracFrenkel} is equivalent to 
 \begin{align}
  \dot{\tenY} = \dot{\tenC} \bigtimes_{n=1}^d \Un + \sum_{n=1}^d \tenC \times_n \dUn \bigtimes_{k \neq n} \Uk
 \end{align}
with the core and the factors satisfying 
\begin{eqnarray}\label{eqn:kochdiff}
\begin{aligned}
  \dot{\tenC} = &\, \dot{\tenA} \bigtimes_{k=1}^d \Ukt, \\
  \dUn   = &\, \Pno \Big[ \dot{\tenA} \bigtimes_{k \neq n} \Ukt \Big]_{(n)}\, \Cnt\big(\Cn \Cnt\big)^{-1}  
\end{aligned}
\end{eqnarray}
where $\Pno = \Id - \Un \Unt$ is the orthogonal projection onto $\RgUno$. 
\end{thm}
The equations \eqref{eqn:kochdiff} can be integrated for a given initial tensor $\tenY(0) \in \Mr$ to obtain the Tucker decompositions $\tenY(t)\in \Mr,\,t>0$, which are approximations to 
$\tenA(t) \in \mathbb{R}^{I_1 \times I_2 \times \cdots \times I_d}$ in the sense of the Galerkin condition \eqref{eqn:DiracFrenkel}, or equivalently, the optimization problem
\begin{align}
 \min_{\dot{\tenY}(t) \in \subTYMr} \|\dot{\tenY}(t) - \dot{\tenA}(t) \|.
\end{align}
The error of \eqref{eqn:DiracFrenkel} can be bounded in terms of the best-approximation error \eqref{eqn:bestappr}, which can be extended to the case of tensor differential equations if (one-sided) Lipschitz conditions for $F$ 
are assumed \cite{koch2010dynamical}.  

\section{Optimization approach and Euler scheme via ALS}\label{sec:alseuler}
\subsection{Dynamical approximation without orthogonality constraint in the tangent space}
We can derive generalizations of \eqref{eqn:kochdiff} if we give up the orthogonality constraint $\Unt\dUn = 0$. We will then numerically tackle these {ungauged} equations by an optimization scheme {for the 
defect approximation \eqref{eqn:optgen} circumventing the explicit form of the factor equations. This will allow us to overcome numerical problems} coming from the orthogonal projections and the necessary inversion of core unfoldings. 
\begin{thm}
 For a tensor $\tenY \in \Mr$ with $n$-mode factors having orthonormal columns the optimization problem
 \begin{align}\label{eqn:optgen}
   \min_{\dot{\tenY}(t) \in \TYMr} \|\dot{\tenY}(t) - \dot{\tenA}(t) \|
 \end{align}
 is equivalent to 
 \begin{align}\label{eqn:tangentten}
  \dot{\tenY} = \dot{\tenC} \bigtimes_{n=1}^d \Un + \sum_{n=1}^d \tenC \times_n \dUn \bigtimes_{k \neq n} \Uk
 \end{align}
with the core and the factors satisfying 
\begin{eqnarray}\label{eqn:diffgen}
\begin{aligned}
  \dot{\tenC} = &\, \dot{\tenA} \bigtimes_{k=1}^d \Ukt - \sum_{k=1}^d \tenC \times_k \Ukt \dUk, \\
  \dUn   = &\, \Big(\Big[ \dot{\tenA} \bigtimes_{k \neq n} \Ukt \Big]_{(n)}\, - \Un \dCn - \Un \sum_{k \neq n} \big[ \tenC \times_k \Ukt \dUk \big]_{(n)} \Big) \Cnt\big(\Cn \Cnt\big)^{-1},
\end{aligned}
\end{eqnarray}
which reduces to \eqref{eqn:kochdiff} if $\Unt\dUn = 0,\, n=1,\hdots,d$. 
\end{thm}
\textit{Proof.}
We look at the optimization problem 
\begin{align*}
 \min_{\dot{\tenY}(t) \in \TYMr} \|\dot{\tenY}(t) - \dot{\tenA}(t) \|.
\end{align*}
Since a tangent tensor in $\TYMr$ has the form \eqref{eqn:tangentten} we will minimize for given $\tenY(t) \in \Mr$ the objective function
\begin{align*}
 \| \big(\dot{\tenC} \bigtimes_{k=1}^d \Uk + \sum_{k=1}^d \tenC \times_k \dUk \bigtimes_{\ell \neq k} \Uell\big) - \dot{\tenA} \|^2
\end{align*}
 with respect to $\dot{\tenC}$ and $\dUk$. Rewriting yields
 \begin{align*}
  \| \dot{\tenC} \bigtimes_{k=1}^d \Uk + \sum_{k=1}^d \tenC \times_k \dUk \bigtimes_{\ell \neq k} \Uell - \dot{\tenA} \|^2 = 
  \|\dot{\tenA}\|^2 + \|\dot{\tenC}\|^2 - 2\langle \dot{\tenA} \bigtimes_{k=1}^d \Ukt, \dot{\tenC} \rangle 
   + \sum_{k=1}^d \|\tenC \times_k \dUk\|^2  &\\
  - 2 \sum_{k=1}^d \langle\dot{\tenA} \bigtimes_{\ell \neq k} \Uellt, \tenC \times_k \dUk \rangle 
  + 2 \sum_{k=1}^d \langle \dot{\tenC}, \tenC \times_k \Ukt \dUk \rangle 
  + \sum_{k=1}^d \sum_{\ell \neq k} \langle \tenC \times_k \Ukt \dUk, \tenC \times_\ell \Uellt \dUell \rangle. &
 \end{align*}
This is a convex quadratic function in $\dot{\tenC}$, hence, setting the gradient with respect to $\dot{\tenC}$ equal to zero gives
the first equation in \eqref{eqn:diffgen}. We now define the part of the objective function involving only terms in $\dUn$ with fixed index $n$ (utilizing $n$th unfolding) as the convex quadratic function
\begin{eqnarray}\label{eqn:fn}
\begin{aligned}
 f_n(\dUn) := \|\dUn \Cn \|^2 - 2\langle \dUn, \big[\dot{\tenA} \bigtimes_{\ell \neq n} \Uellt]_{(n)} \Cnt \rangle
 + 2 \langle \dUn, \Un \dCn \Cnt \rangle \\+ 2 \sum_{\ell \neq n} \langle \dUn, \Un \big[\tenC \times_\ell \Uellt \dUell\big]_{(n)} \Cnt \rangle .
\end{aligned}
\end{eqnarray}
The derivative of $f_n$ set to zero gives the second equation of \eqref{eqn:diffgen}. \hfill $\Square$\\

Note that a discretization scheme using the derivatives \eqref{eqn:kochdiff} or \eqref{eqn:diffgen} requires the solution of a least squares problem 
for the factor matrices by calculating the pseudo inverse $\Cnt\big(\Cn \Cnt\big)^{-1}$ of the $n$th-unfolding of the core tensor. 
Small singular values in $\Cn$ will affect the solution and make the method unstable \cite{kieri2016discretized}. However, the advantage of \eqref{eqn:diffgen} lies in the fact that 
the derivative $\dUn$ is not orthogonal to the range of $\Un$. This becomes apparent in a discretization scheme that we will establish in the forthcoming {section, where we will establish a numerical time-stepping scheme that, 
in particular, circumvents the explicit formulation of the factor differential equations of \eqref{eqn:diffgen}. }
\subsection{Euler method under {orthogonal columns} constraint}
We will consider Euler's method for {solving \eqref{eqn:optgen}} via an alternating least squares (ALS) scheme. 
For that purpose, all but one factor matrix derivative will be kept fixed consecutively, where the core derivative is assumed to be already computed according to the first equation in \eqref{eqn:diffgen}.  
Consider an Euler step at $t$ for the second equation of \eqref{eqn:diffgen} with length $h>0$, that is  
\begin{align}\label{eqn:Euler}
 \Un(t+h) := \Un(t) + h\, \dUn(t).  
\end{align}
Note that the orthogonality gauge $\Unt\dUn = 0$ prevents the 
updated factor matrices from having orthonormal columns, that would be $\Un(t+h)^T \Un(t+h) = \Id$. 
On the other hand, as will become apparent from Theorem~\ref{thm:trace} below, the function $f_n$ \eqref{eqn:fn} with respect to the variable $\Unh :=\Un(t+h)$ gets linear 
if ${\Unh}^T \Unh = \Id$ is imposed. 
We will therefore seek for an Euler step under the orthogonality constraint for the updated factor matrices. 
More precisely, for fixed $t$ we define for  $\Vn \in \mathbb{R}^{I_n \times r_n}$ (omitting the argument $t$)
\begin{align}
 \widetilde{f}_n(\Vn):=f_n\big(\tfrac{1}{h}(\Vn-\Un)\big),
\end{align}
and define $\Unh$ as the solution of the constraint optimization problem 
\begin{align}\label{eqn:opteuler}
 \min_{\Vnt \Vn = \Id} \widetilde{f}_n(\Vn),
\end{align}
{where the involved derivative $\dot{\tenA}(t)$ is treated explicitly.}  
The following theorem characterizes the updated factor matrices $\Unh$ from one Euler step according to \eqref{eqn:opteuler} as the solution to a trace optimization problem, 
which is computable from an SVD of a  matrix of size $I_n \times r_n$.
\begin{thm}\label{thm:trace}
 The optimization problem
 \begin{align}\label{eqn:minfn}
  \min_{\Vnt\Vn = \Id} \widetilde{f}_n(\Vn)
 \end{align}
 is equivalent to the trace optimization problem
 \begin{align}\label{eqn:traceopt}
  \max_{\Vnt\Vn = \Id} \mathrm{Tr}\big[\Bnt \Vn\big],  
 \end{align}
with 
\begin{align}\label{eqn:Bn}
 \Bn =  \Big(\Big[ \dot{\tenA} \bigtimes_{k \neq n} \Ukt \Big]_{(n)}\, - \Un \dCn - \Un \sum_{k \neq n} \big[ \tenC \times_k \Ukt \dUk \big]_{(n)} \Big) \Cnt + \tfrac{1}{h} \Un \Cn \Cnt.
\end{align}
Moreover, given the (economy size) SVD $\Bn = \Rn \mathbf{\Sigma}_n \Tnt$ with $\Rn \in \mathbb{R}^{I_n \times r_n}, \Tn \in \mathbb{R}^{r_n \times r_n}$ and the diagonal matrix $\mathbf{\Sigma}_n \in \mathbb{R}^{r_n \times r_n}$, 
the solution $\Unh$ to \eqref{eqn:traceopt} is given by
\begin{align}
 \mathbf{V}^\ast_n \equiv \Unh = \Rn \Tnt.  
\end{align}
\end{thm}
\textit{Proof.}
Let $\Bnp = \Big(\Big[ \dot{\tenA} \bigtimes_{k \neq n} \Ukt \Big]_{(n)}\, - \Un \dCn - \Un \sum_{k \neq n} \big[ \tenC \times_k \Ukt \dUk \big]_{(n)} \Big) \Cnt$. Then we have 
$f_n(\dUn) = \|\dUn \Cn\|^2 -  2\langle \dUn, \Bnp\rangle$ and hence
\begin{align*}
 \widetilde{f}_n(\Vn) = \tfrac{2}{h^2}\big( \|\Cn\|^2 - \langle \Vn, \Un \Cn \Cnt\rangle \big) - \tfrac{2}{h} \big( \langle \Vn,\Bnp\rangle - \langle \Un,\Bnp\rangle\big),
\end{align*}
where we used $\Vnt\Vn = \Unt\Un = \Id$. Thus, the optimization problem \eqref{eqn:minfn} is equivalent to 
\begin{align*}
 \max_{\Vnt\Vn = \Id} \langle \Vn,\Bnp + \tfrac{1}{h} \Un \Cn \Cnt\rangle,
\end{align*}
which is the trace optimization problem \eqref{eqn:traceopt}. Consider now the singular value decomposition $\Bn = \tRn \widetilde{\mathbf{\Sigma}}_n \Tnt$ with 
$\tRn \in \mathbb{R}^{I_n \times I_n}$,  
$\widetilde{\mathbf{\Sigma}}_n = 
\left[\begin{array}{c}
  \mathbf{\Sigma}_n\\
  \mathbf{0}
\end{array}\right] = 
\left[\begin{array}{c}
  \mathrm{diag}(\sigma_1, \hdots,\sigma_{r_n})\\
  \mathbf{0}
\end{array}\right] \in \mathbb{R}^{I_n \times r_n}$ and $\Tn \in \mathbb{R}^{r_n \times r_n}$. Note that we assume $r_n \leq I_n$. 
Cyclic permutation of the arguments in the trace yields
\begin{align*}
  \mathrm{Tr}\big[\Bnt \Vn\big] =  \mathrm{Tr}\big[\Tn \widetilde{\mathbf{\Sigma}}_n \tRnt \Vn\big] =  \mathrm{Tr}\big[\widetilde{\mathbf{\Sigma}}_n \tRnt \Vn\Tn \big] = \sum_{j=1}^{r_n} \sigma_j w_{jj},
\end{align*}
where we define the matrix $\Wn := \tRnt \Vn\Tn \in \mathbb{R}^{I_n \times r_n}$. Owing to $\Wnt \Wn = \Id$, we have for the diagonal elements $|w_{jj}|\leq 1$. 
Hence, the maximum of the objective in \eqref{eqn:traceopt} is attained for $w_{jj} \equiv 1$, which implies $\Id = \Rnt \Vn \Tn$ and consequently $\Vn = \Rn \Tnt $. 
\hfill $\Square$\\
 
Note that the matrices $\Bn$ do not require inversion of the potentially problematic matrices $\Cn \Cnt$.\\
Theorem~\ref{thm:trace} readily generalizes to the case of complex valued dynamical tensor approximation. More precisely, 
in the subproblems of the ALS scheme the functions $f_n$ boil down to the real part of the inner product, that is $\mathrm{Re}[\langle \dUn , \Bnp \rangle]$. As a consequence, the equivalent 
trace optimization problem takes the form
\begin{align}
 \max_{{\Ab}^\ast \Ab = \Id} \mathrm{Re}[\langle \Ab , \Bb \rangle],
\end{align}   
with the solution 
\begin{align}
 \Ab_0 = \Rb \Tb^{\text{H}},
\end{align}
where $\Bb = \Rb \mathbf{\Sigma} \Tb^{\text{H}}$ the economy sized SVD of $\Bb$. 

\subsection{Algorithm for Euler step via ALS}
The computation of an Euler step for \eqref{eqn:diffgen} via alternating least squares for \eqref{eqn:optgen} is summarized in Algorithm~\ref{alg:alseuler}.

\begin{algorithm}\caption{Euler step from $t \rightarrow t+h$ via ALS for \eqref{eqn:optgen}}\label{alg:alseuler}
 {Given $\tenY(t) = \tenC(t) \bigtimes_{n=1}^d \Un(t),\, \dot{\tenA} = \dot{\tenA}(t)$.}\\
 Initialize $\deltaUn$, e.g., $\deltaUn = 0,\, n=1\hdots,d$.\\
 Initialize $\Delta{\tenC} = \dot{\tenA} \bigtimes_{k=1}^d \Ukt - \sum_{k=1}^d \tenC \times_k \Ukt \deltaUk$.\\
 \textbf{repeat}
 \begin{algorithmic}
 \For{$n= 1,\hdots,d$}
 \State Compute $\Bn$ from \eqref{eqn:Bn} with $\dUn \gets \deltaUn$ and $\dot{\tenC} \gets \Delta{\tenC}$
 \State Compute economy SVD $\Bn = \Rn \mathbf{\Sigma}_n \Tnt$
 \State Set $\Unh = \Rn \Tnt$
 \State Set $\deltaUn = (\Unh -\Un)/h$
 \EndFor
 \State Update core  $\Delta{\tenC} = \dot{\tenA} \bigtimes_{k=1}^d \Ukt - \sum_{k=1}^d \tenC \times_k \Ukt \deltaUk$
 \end{algorithmic} 
 \textbf{until} fit \eqref{eqn:fit} ceases to improve or maximum iterations exhausted\\
 \textbf{return} Approximation to $\tenY(t+h) \approx (\tenC + h\Delta{\tenC}) \bigtimes_{n=1}^d \Unh$
\end{algorithm}

An Euler method for the dynamical approximation \eqref{eqn:optgen} using Algorithm~\ref{alg:alseuler} starts from a rank-$\mathbf{r}$ decomposition or (best-) approximation of $\tenA(0)$, e.g., 
via higher order orthogonal iteration (HOOI), see \cite{kolda2009tensor}. 
{The termination criterion uses the discrete version of the norm of the defect (resp. residual) \eqref{eqn:optgen}, that is
\begin{align}\label{eqn:fit}
  \| \big(\Delta{\tenC} \bigtimes_{k=1}^d \Uk + \sum_{k=1}^d \tenC \times_k \deltaUk \bigtimes_{\ell \neq k} \Uell\big) - \dot{\tenA} \|.
\end{align}}
{Algorithm~\ref{alg:alseuler} computes one time step of the underlying \textit{explicit Euler scheme} for the initial value problem associated with the defect approximation \eqref{eqn:optgen}.} 
Higher order schemes of Runge-Kutta type can be composed based on the above Euler scheme, {see section~\ref{sec:secord}.} 
{The derivative at time $t$, i.e., $\dot{\tenA}(t)$ or $F(t,\tenY(t))$ respectively, is treated explicitly and assumed to be efficiently computable. 
The overall computational effort depends heavily on this derivative evaluation, compare for instance with the numerical example in section~\ref{sec:rankdef}. One inner loop evaluation in the Euler-ALS scheme described 
in Algorithm~\ref{alg:alseuler} involves the computation of the rectangular matrix $\Bn$ of size $I_n \times r_n$ at the costs of $\mathcal{O}\big((\prod_{k=1}^d r_k) \,\sum_{k=1}^d I_k\big)$ operations plus 
the costs for the term $[ \dot{\tenA} \bigtimes_{k \neq n} \Ukt ]_{(n)}$. If we assume $\dot{\tenA} \in \mathscr{M}_{\mathbf{r}^\prime}$ (e.g. after projection/compression), the latter term takes 
$\mathcal{O}\big(I_n r_n^\prime \prod_{k \neq n} r_k + (\prod_{k \neq n} r_k^\prime) \,\sum_{k \neq n} r_k\big)$ operations. The economy sized SVD scales with $\mathcal{O}(I_n r_n^2)$.}

\subsection{Regularization}
In the course of the derivation of Algorithm~\ref{alg:alseuler} the orthogonality gauge $\Unt \dUn = 0$ of \cite{koch2010dynamical} was given up. 
As a consequence, elements of the tangent space $\TYMr$ have no unique representation any more. Although, we have not observed any related problems in the numerical tests, we found that including a regularization of 
Tikhonov type helps to stabilize the SVD in cases where $\Bn$ is effectively rank deficient and where $\Cn \Cnt$ is numerically close to singular, respectively. More precisely, we can treat the problem ($\alpha > 0$)
\begin{align}
 \min_\dUn f_n(\dUn) + \alpha \, \|\dUn\|^2,
\end{align}
compare with \eqref{eqn:fn}. Note that this positively affects definiteness of the quadratic term in \eqref{eqn:fn}, which becomes $\langle\dUn, \dUn (\Cn \Cnt + \alpha \Id) \rangle$.
The discrete scheme changes as follows. We treat the modified optimization problem 
\begin{align}\label{eqn:minfnmod}
 \min_{\Vnt\Vn = \Id} \widetilde{f}_n(\Vn) - \tfrac{2\alpha}{h^2} \langle \Vn, \Un \rangle,
\end{align}
which is equivalent to
\begin{align}\label{eqn:traceoptmod}
 \max_{\Vnt\Vn = \Id} \mathrm{Tr}\big[ (\Bn + \tfrac{\alpha}{h}\, \Un )^T\Vn\big]  
\end{align}
according to Theorem~\ref{thm:trace}.
Observe that the regularization term is only effective because $\Unt \dUn \neq 0$. In fact, if $\dUn = \tfrac{1}{h}(\Vn - \Un) \in \RgUno$ we get 
\begin{align*}
 \Pn \Vn = \Un,
\end{align*}
where $\Pn$ is the orthogonal projection onto $\RgUn$. Hence, the regularization term would degenerate to 
$\mathrm{Tr}[\Unt \Vn] =  \mathrm{Tr}[\Unt (\Pn\Vn + \Pno\Vn)] =  \mathrm{Tr}[\Unt \Pn\Vn] = \|\Un\|^2$.\\ 

The choice $\alpha\, \propto\, h^2$ makes the regularization term in \eqref{eqn:minfnmod} independent of the step size $h$. If regularization was considered in the numerical examples below we have always chosen $\alpha = h^2$.

\section{Numerical validation}
We use the Tensor Toolbox Version 2.6 \cite{TTB_Software}, which implements efficient tensor arithmetics and fitting algorithms for different tensor classes in MATLAB \cite{TTB_Sparse}.
Results were achieved by using the Vienna Scientific Cluster 3 (VSC3). 

\subsection{Koch/Lubich Example}
\begin{figure}[hbtp]
\hspace{0.0 in} 
\centering 
\includegraphics[scale=0.35]{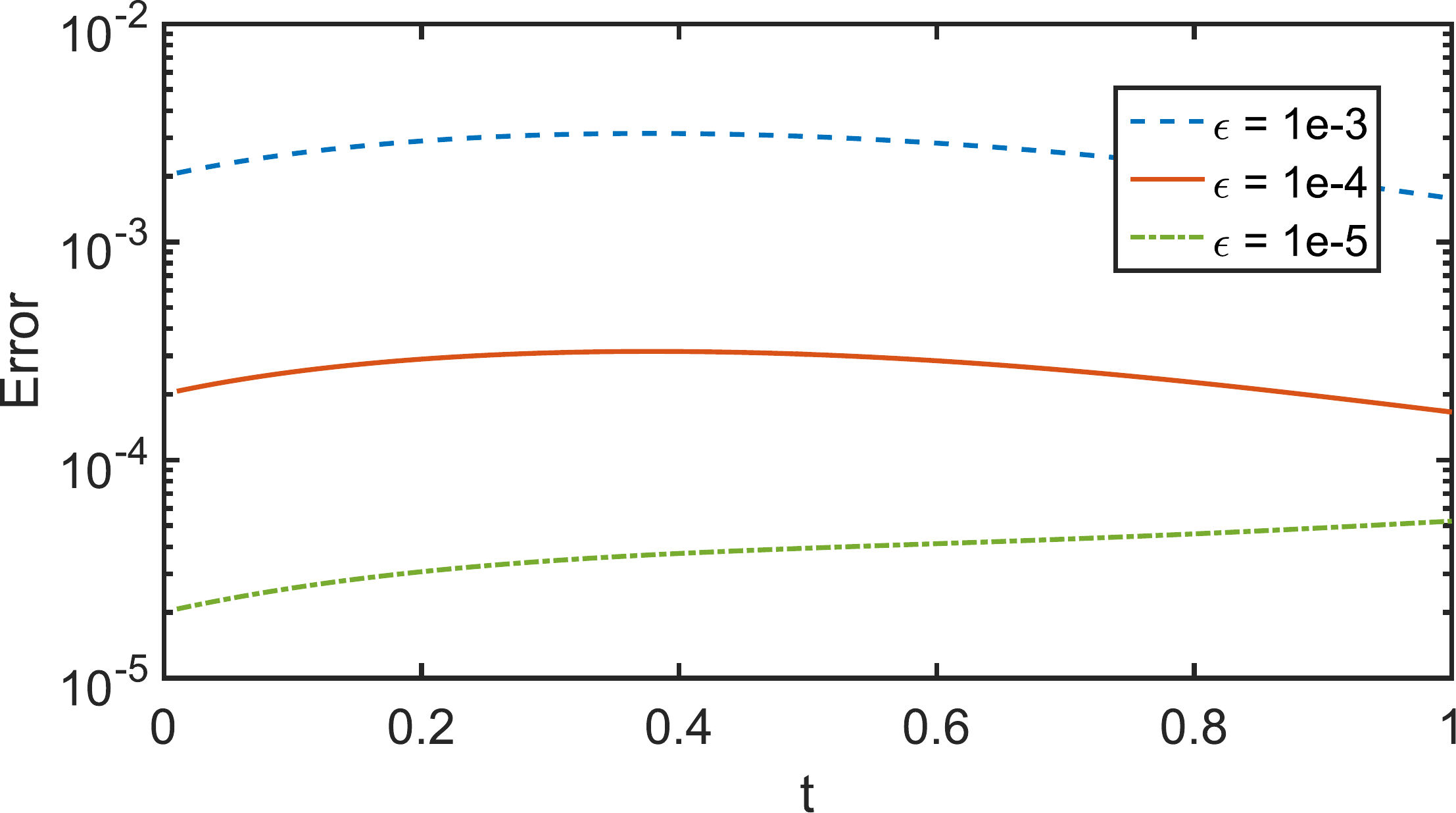}
\caption{Errors of the dynamical approximation of \eqref{eqn:expl1} with $\mathbf{r} = (10,10,10,10)^T$ with $\varepsilon \in \{ 10^{-5}, 10^{-4}, 10^{-3}\}$. Interval $[0,1]$, step size $h=10^{-4}$.}\label{fig:kochlubich1}
\end{figure}

We consider the numerical example as in \cite{koch2010dynamical}.  
A time-dependent data tensor $\tenA \in \mathbb{R}^{15 \times 15 \times 15 \times 15}$ is constructed as 
\begin{align}\label{eqn:expl1}
 \tenA(t) = \mathrm{exp}(t) \, \tenB \bigtimes_{k=1}^4 \Uk + \varepsilon \, \big(1+t+\sin(3\,t)\big)\,\tenC, \quad t \in [0,1],
\end{align}
where $\Uk \in \mathbb{R}^{15 \times 10}$ are random matrices with orthonormal columns, $\tenB \in \mathbb{R}^{10 \times 10 \times 10 \times 10}$ a random core tensor, and $\tenC \in \mathbb{R}^{15 \times 15 \times 15 \times 15}$ 
a random perturbation. The dynamical Tucker tensor approximations are computed with $\mathbf{r} = (10,10,10,10)^T$, fixed $h=10^{-4}$ and different $\varepsilon \in \{ 10^{-5}, 10^{-4}, 10^{-3}\}$. 
We use Euler steps described in Algorithm~\ref{alg:alseuler} without regularization. For initialization we use $\deltaUn = 0,\, n=1\hdots,d$. Moreover, a quasi best-approximation for $\tenA(0) \approx \tenA_0$ is used, which is derived 
from HOOI. Note that this example solves the tensor initial value problem 
\begin{align*}
 \dot{\tenA} = &\, \tenA + \varepsilon \big( 3\cos(3\,t) - \sin(3\,t) - t\big) \tenC, \\
 \tenA(0) =& \, \tenA_0.
\end{align*}
Results are shown in Figure~\ref{fig:kochlubich1}, where relative errors are given, i.e., $\|\tenA(t) - \tenY(t)\|/\|\tenA(t)\|$. Alg.~\ref{alg:alseuler} was used with a tolerance of $10^{-5}$ for the change in the relative {defect} 
$\|\dot{\tenA}(t) - \dot{\tenY}(t)\|/\|\dot{\tenA}(t)\|$ (cf. \eqref{eqn:fit}), which was basically reached after two loops. {The norm of the relative defect itself was constantly in the order of $\varepsilon$ 
(about $5 \varepsilon$).} 

\subsection{Second Order Method}\label{sec:secord}
\begin{figure}[hbtp]
\hspace{0.0 in} 
\centering 
\includegraphics[scale=0.35]{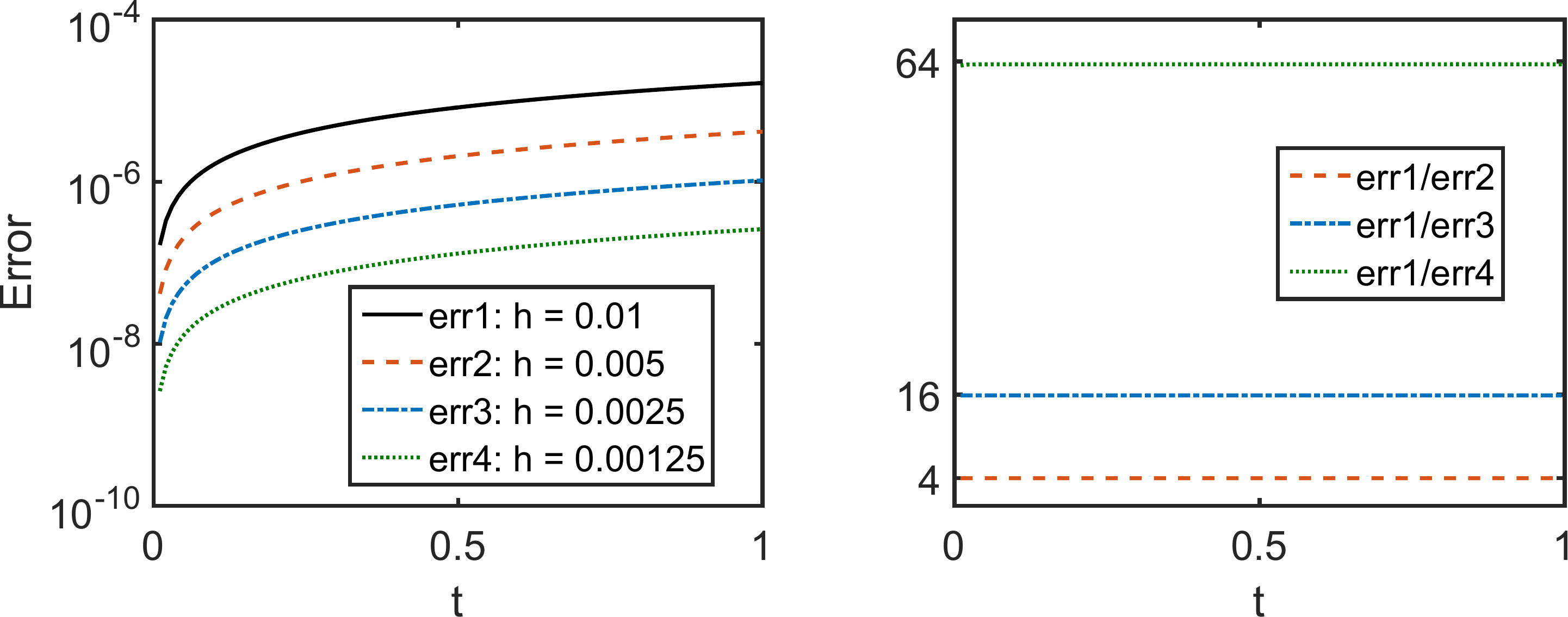}
\caption{Errors of the dynamical approximation of \eqref{eqn:expl1} with $\mathbf{r} = (8,9,10,11)^T$ and $\varepsilon = 10^{-10}$ using the improved Euler method (left) and comparison of the errors among each other (right).}\label{fig:kochlubich_impr}
\end{figure}

We repeat Example 1 with $\mathbf{r} = (8,9,10,11)^T$ and adjusted sizes in \eqref{eqn:expl1} for $\Uk$ and $\tenB$. The second order improved Euler method (extrapolated Euler) is now used, i.e., for an initial value problem
\begin{align*}
 y^\prime &\,= f(t,y)\\
 y(0) &\, = \eta_0
\end{align*}
the new iterate $\eta_\nu$ at time $t_\nu = t_{\nu-1} + h$ is calculated from the previous iterate $\eta_{\nu-1}$ at time $t_{\nu-1}$ according to
\begin{align*}
 \eta_\nu = \eta_{\nu-1} +  h f(t_{\nu-1} + \tfrac{h}{2}, \eta_{\nu-1} + \tfrac{h}{2} f(t_{\nu-1},\eta_{\nu-1}) ),
\end{align*}
which is a $2$-stage explicit Runge-Kutta formula. The half and full Euler steps, respectively, are calculated using Alg.~\ref{alg:alseuler}. 
Figure~\ref{fig:kochlubich_impr} shows relative errors for different step sizes $h$ and $\varepsilon = 10^{-10}$ on the interval $[0,1]$ (left) and 
compares these errors among each other (right).
We can observe second order convergence in $h$. {Again the norm of the relative defect for the time steps was in the order of about $\varepsilon$ independent of the step size.}

\subsection{Small singular values in higher dimensions}\label{sec:smallsv}
\begin{figure}[hbtp]
\hspace{0.0 in} 
\centering 
\includegraphics[scale=0.35]{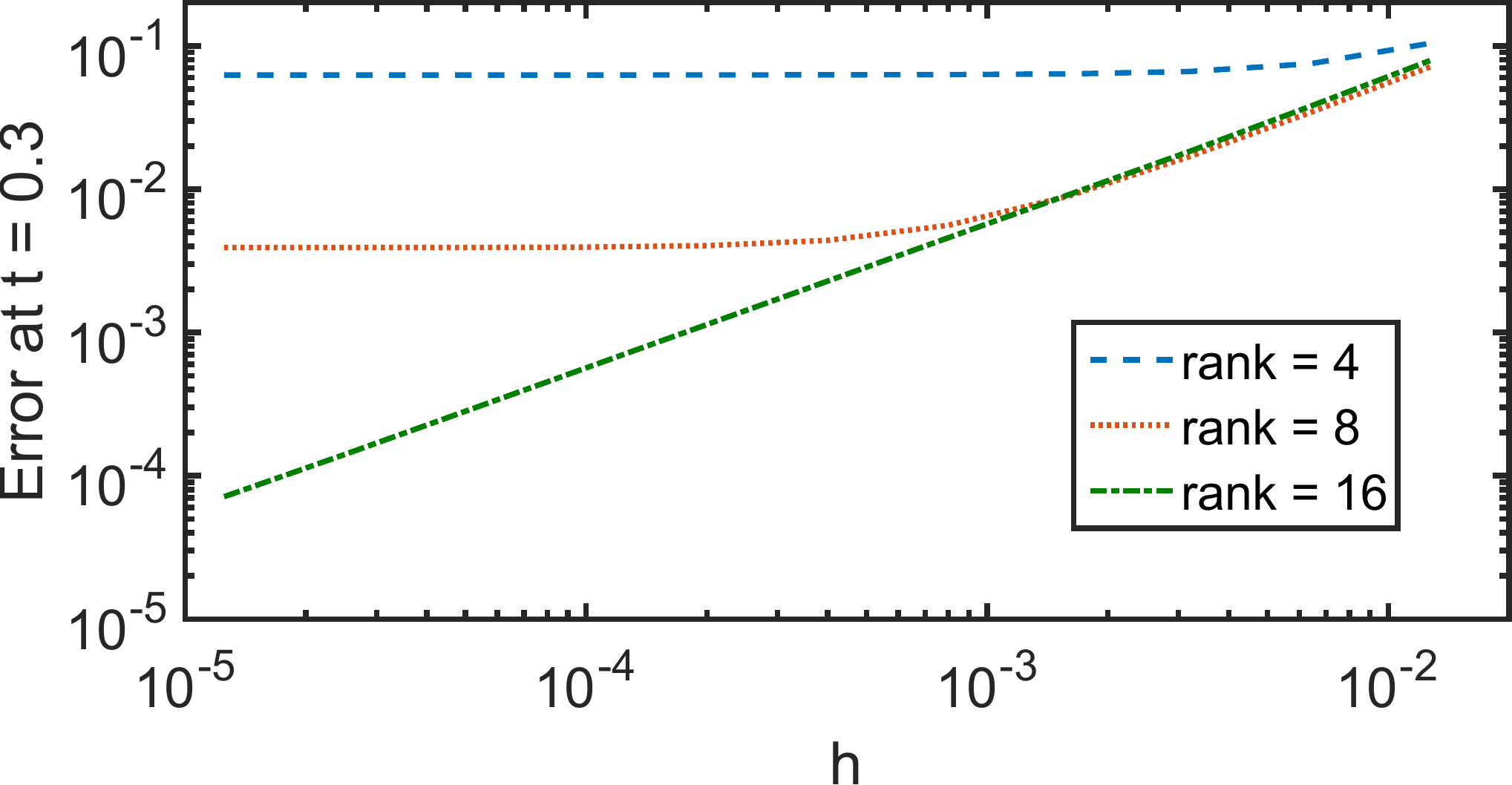}
\caption{Errors at $t=0.3$ of the dynamical approximation of \eqref{eqn:lubichho} in the case $d=2$ and $I=100$ for different step sizes $h$ and ranks $\mathbf{r} = (r,r)$.}\label{fig:lubichho2d}
\end{figure}

We generalize the example from \cite{kieri2016discretized} to higher dimensions. A time-dependent tensor $\tenA(t) \in \mathbb{R}^{I_1 \times I_2 \times \cdots \times I_d}$ with $I_k \equiv I$ is constructed as
\begin{align}\label{eqn:lubichho}
 \tenA(t) =  \mathrm{exp}(t) \, \tenC \bigtimes_{k=1}^d \mathrm{exp}(t\mathbf{W}_k)
\end{align}
with skew-symmetric random matrices $\mathbf{W}_k \in \mathbb{R}^{I \times I}$ and a super-diagonal tensor $\tenC \in \mathbb{R}^{I_1 \times I_2 \times \cdots \times I_d}$ with 
diagonal entries $c_j = 1/2^{(d-1)j},\, j=1,\hdots,I$. We now use the Euler method of Alg.~\ref{alg:alseuler} and regularization with $\alpha = h^2$. Figure~\ref{fig:lubichho2d} shows relative errors at $t=0.3$ against step length $h$ for different ranks in the two dimensional case 
with $I = 100$. 

\begin{figure}[hbtp]
\hspace{0.0 in} 
\centering 
\includegraphics[scale=0.35]{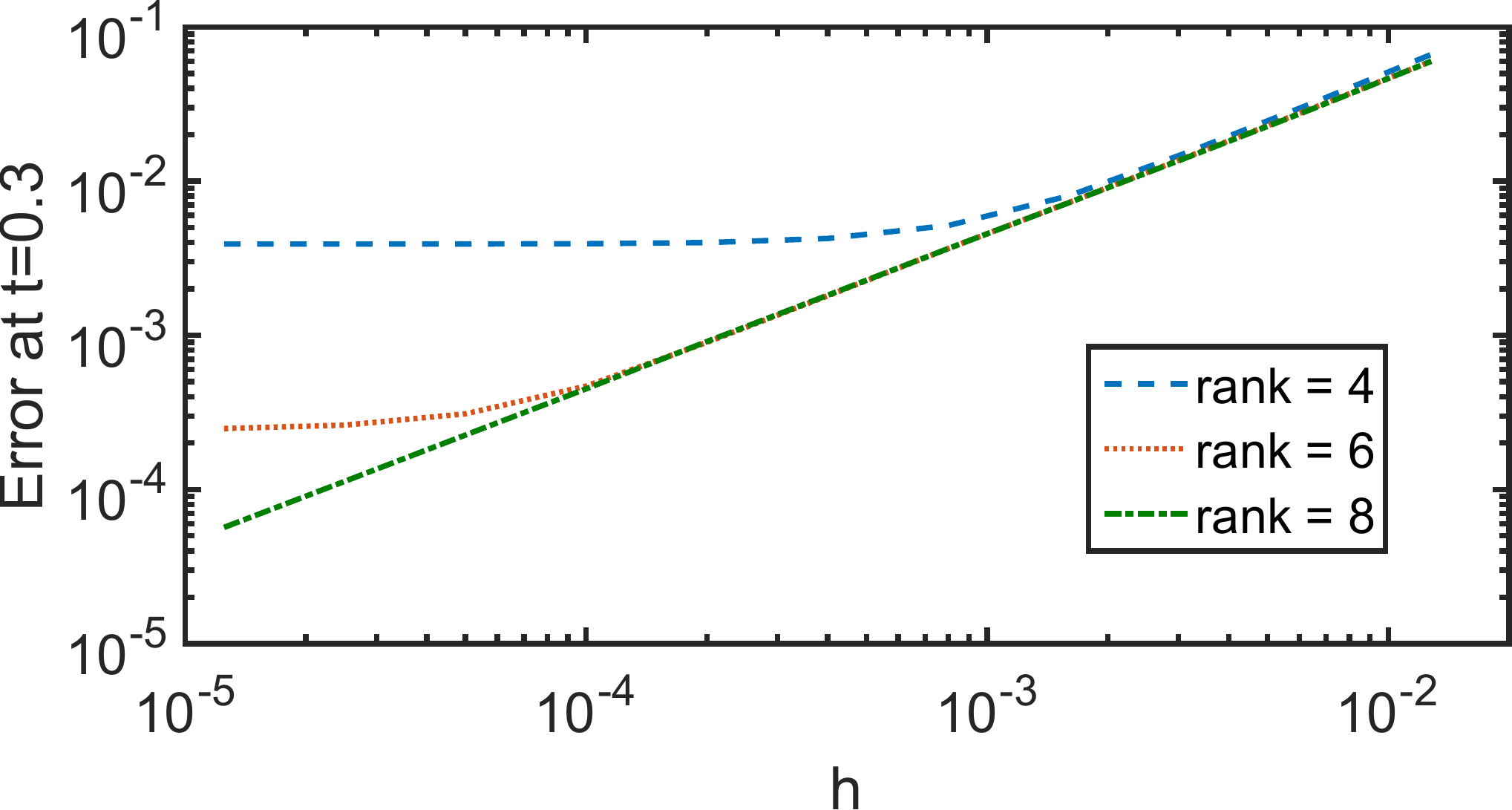}
\caption{Errors at $t=0.3$ of the dynamical approximation of \eqref{eqn:lubichho} in the case $d=3$ and $I=50$ for different step sizes $h$ and ranks $\mathbf{r} = (r,r,r)$.}\label{fig:lubichho3d}
\end{figure}

We repeat these computations for $d=3$ and $I=50$ and give associated results in Figure~\ref{fig:lubichho3d}.

The tolerance for the change of the fit was $10^{-6}$ and basically reached after two loops in Algorithm~\ref{alg:alseuler}.\\ 
{We observe in both cases the desired behavior of linearly decreasing error with respect to $h$ towards the model accuracy, and stability also for higher ranks and larger step size.}

{
\subsection{Spatially discretized PDEs}\label{sec:rankdef}
We consider the heat equation in $2$d with a source term similar as in \cite{Kieri_V:2017}:
\begin{eqnarray}\label{eqn:heateqn2d}
\begin{aligned}
 \frac{\partial u}{\partial t} &\,= 0.01 \,\Delta u + x^{(1)}x^{(2)}\,e^{-t\,(x^{(1)}+x^{(2)})}  \in \Omega = (0,1)^2, \,t > 0, \\
  u(\mathbf{x},t) &\, = 0, \quad \mathbf{x} \in \partial \Omega,\, t>0,\\
  u(\mathbf{x},0) &\,=  e^{-100\,\big((x^{(1)}-0.5)^2+(x^{(2)}-0.5)^2\big)},\, \mathbf{x}\in \Omega.
\end{aligned}
\end{eqnarray}
For numerical computation we use the forward in time and central in space (FTCS) scheme. Uniform spatial discretization of the Laplace operator by second order centered finite differences leading to a \textit{Kronecker sum} 
\begin{align}
 \Delta \approx \mathbf{K}_1 \otimes \Id \,+\, \Id \otimes \mathbf{K}_1,
\end{align}
with $ \mathbf{K}_1,\, \Id\in \mathbb{R}^{I \times I}$ and $K_1 = 1/k^2\,\text{tridiag}(1,-2,1),\,k = 1/I$. Multiplication of such a discretization with a low-rank grid function can be performed efficiently. 
In our numerical experiments we use $I = 200$ and Alg.~\ref{alg:alseuler} with regularization for forward time stepping. 
The source term in \eqref{eqn:heateqn2d} leads to a rapidly increasing effective rank of the numerical solution. 
\begin{figure}[hbtp]
\hspace{0.0 in} 
\centering 
\includegraphics[scale=0.35]{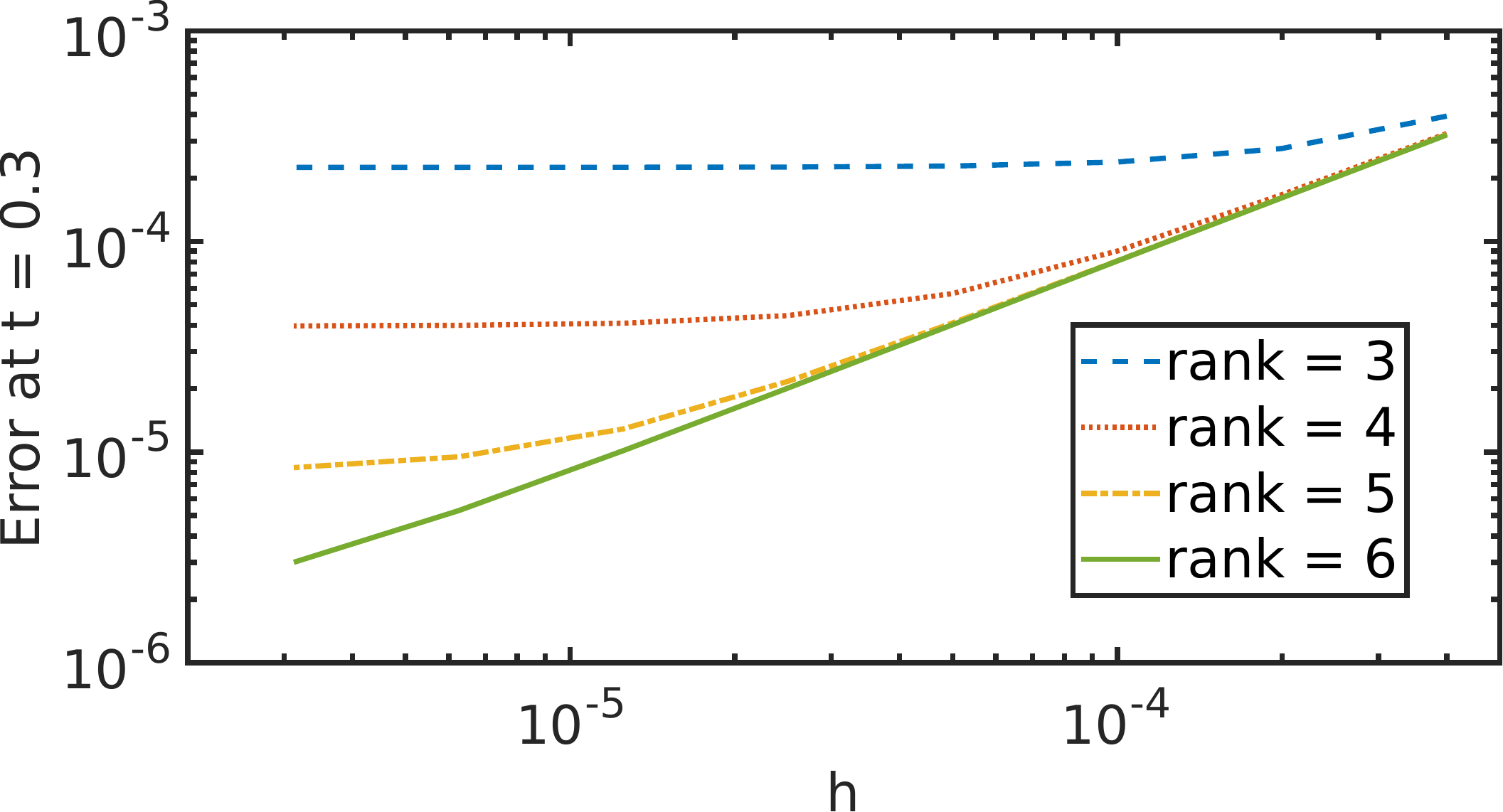}
\caption{Errors at $t=0.3$ of the dynamical approximation of \eqref{eqn:heateqn2d} in the case $d=2$ and $I=200$ for different step sizes $h$ and ranks $\mathbf{r} = (r,r)$.}\label{fig:heateqn2d}
\end{figure}
Fig.~\ref{fig:heateqn2d} shows for different ranks the errors at time $t=0.3$ with respect to a reference solution 
obtained by fourth order Runge-Kutta for the full matrix equation. Note that stability for the FTCS scheme requires a step size below $h = 6.25\cdot 10^{-4}$. 
We observe a similar behavior as in the previous example of section~\ref{sec:smallsv}.\\  
Next we take a reaction-diffusion example in $d=2,3$ space dimensions similar to \cite{nonnenmacher2008dynamical}:
\begin{eqnarray}\label{eqn:reacdiff}
\begin{aligned}
 \frac{\partial u}{\partial t} &\,= 0.01 \,\Delta u +  0.1\,u^3,\quad \mathbf{x} \in \Omega = (0,1)^d, \,t > 0, \\
  u(\mathbf{x},t) &\, = 0, \quad \mathbf{x} \in \partial \Omega,\, t>0,\\
  u(\mathbf{x},0) &\,= 10d\, \prod_{n=1}^d e^{-100\,(x^{(n)}-0.5)^2},\, \mathbf{x}\in \Omega.
\end{aligned}
\end{eqnarray}
We discretize the Laplace operator uniformly by second order centered finite differences leading to a \textit{$d$-term Kronecker sum} 
\begin{align}
 \Delta \approx \mathbf{K}_1 \otimes \Id \otimes \cdots \otimes \Id \,+\, \Id \otimes \mathbf{K}_1 \otimes \cdots \otimes \Id \,+\, \cdots \,+\, \Id \otimes  \Id \otimes \cdots \otimes  \mathbf{K}_1,
\end{align}
with $ \mathbf{K}_1,\, \Id\in \mathbb{R}^{I \times I}$ and $K_1 = 1/k^2\,\text{tridiag}(1,-2,1),\,k = 1/I$. 
Note that such a discretization of an operator as a sum of Kronecker product matrices applied to a grid function represented in Tucker format results in a sum of Tucker tensors of same rank, e.g., \cite{exl2014fft}. 
This makes computations on large grids possible. The nonlinearity in \eqref{eqn:reacdiff} leads to elementwise (\textit{Hadamard}) products of tensors, which results in cubed Tucker ranks. 
The Hadamard products are recompressed via HOOI, however, more sophisticated efficient techniques for Hadamard products in the Tucker format were recently introduced \cite{kressner2017recompression} but not considered here. 
We solve \eqref{eqn:reacdiff} subject to the low-rank manifold $\Mr$ for different ranks $\mathbf{r} \equiv r \in \{ 3,4,5\}$, where we use 
the improved Euler method of example~\ref{sec:secord} with basic time steps realized by Alg.~\ref{alg:alseuler} using regularization. 
The computations of Figure~\ref{fig:reacdiff2d} were carried out for the case $d=2$ with $I = 400$ and step size $h = 10^{-4}$, where the norm of the relative defect as a measure for local (time-stepping and projection) errors for different ranks is shown during time integration towards the blow-up. 
\begin{figure}[hbtp]
\hspace{0.0 in} 
\centering 
\includegraphics[scale=0.35]{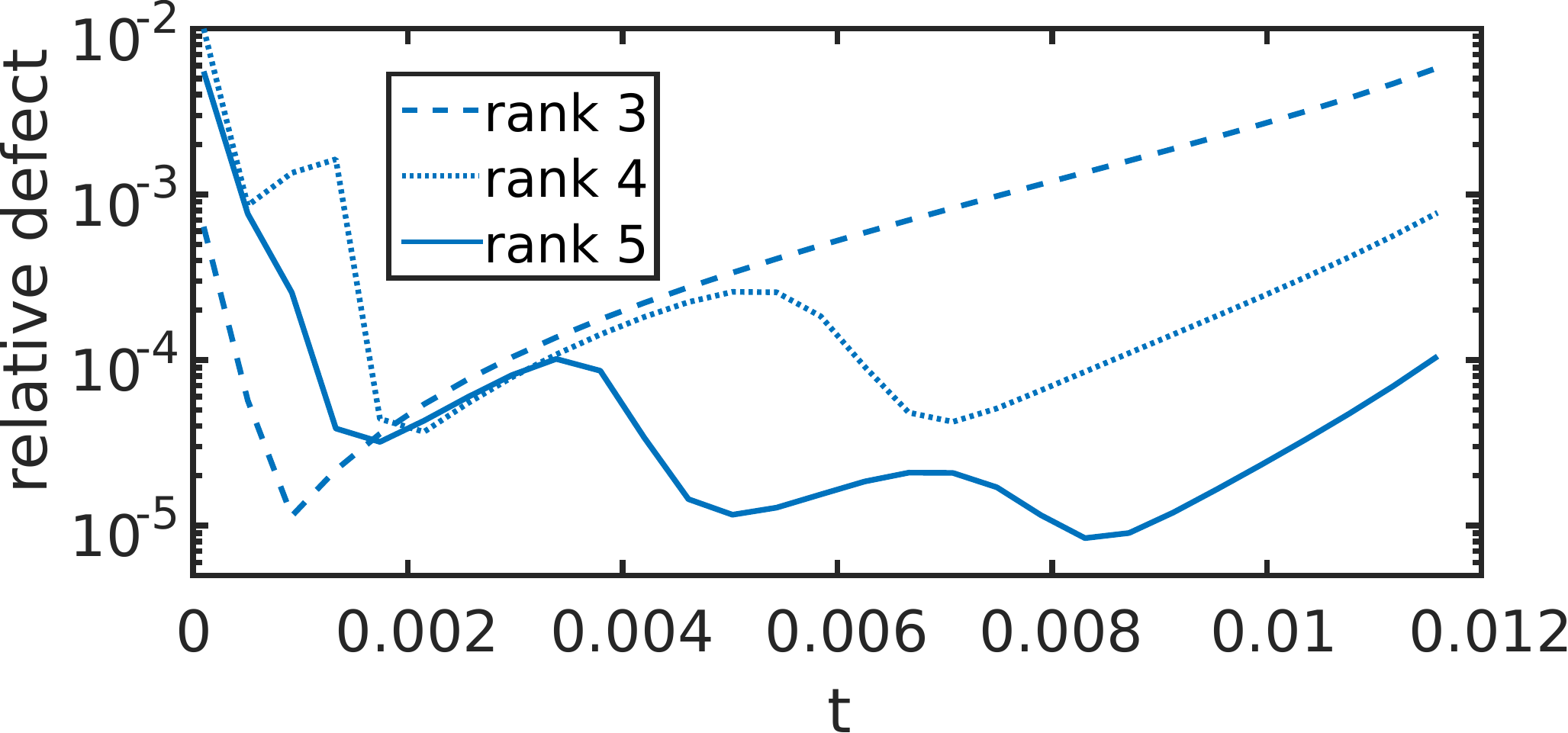}
\caption{Local relative defect of the dynamical approximation of \eqref{eqn:reacdiff} in the case $d=2$ and $I=400$ for different ranks $\mathbf{r} = (r,r)$.}\label{fig:reacdiff2d}
\end{figure}
One observes decreasing local defect in the initial stabilization phase and linearly increasing defect afterwards. The latter is due to 
effective rank growth caused by the nonlinear term in \eqref{eqn:reacdiff}. 
Fig.\ref{fig:blowup} shows the initial data and the numerical solution at time $t=0.012$ corresponding to problem 
\eqref{eqn:reacdiff}. 
\begin{figure}[hbtp]
\hspace{0.0 in} 
\centering 
\includegraphics[scale=0.38]{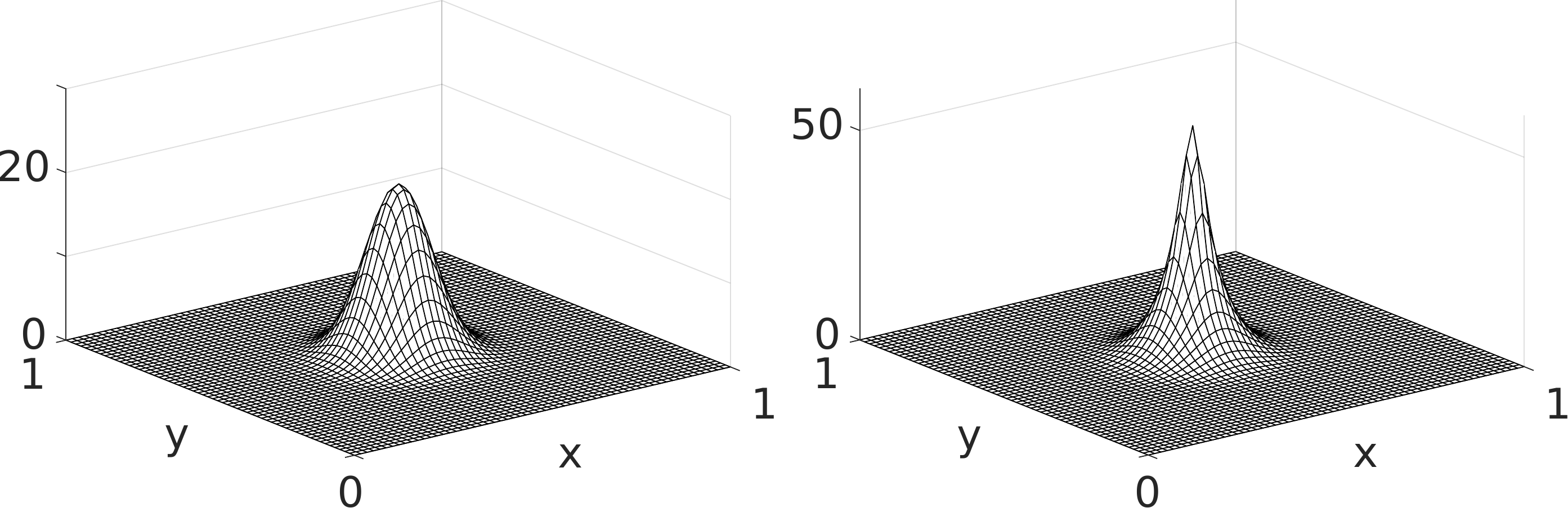}
\caption{Initial data (left) and the numerical solution at time $t=0.012$ (right) corresponding to \eqref{eqn:reacdiff} in the case $d=2$.}\label{fig:blowup}
\end{figure}
Fig.~\ref{fig:compare} shows the difference between a reference solution obtained from an adaptive fourth order Runge-Kutta scheme applied to the full $2$d system and the rank $3$ and $5$ low-rank numerical solutions. 
\begin{figure}[hbtp]
\hspace{0.0 in} 
\centering 
\includegraphics[scale=0.3]{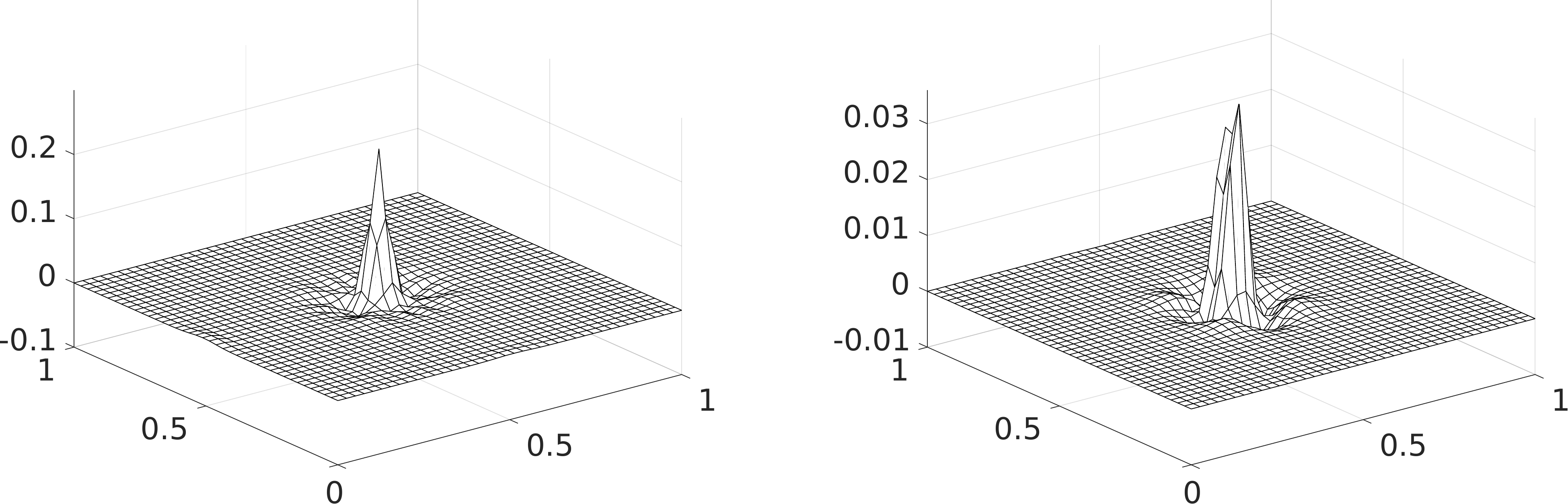}
\caption{Deviation from reference solution for $2d$ system \eqref{eqn:reacdiff} at time $t=0.012$ for dynamical rank-$3$ (left) and rank-$5$ (right) approximation.}\label{fig:compare}
\end{figure}
We repeat the computations for the case of $d=3$, where treatment of the full equations gets to expensive both in storage and computational effort. Fig.~\ref{fig:reacdiff3d} shows local relative defects for different ranks, where we can observe a similar behavior as for the $2$d case. 
\begin{figure}[hbtp]
\hspace{0.0 in} 
\centering 
\includegraphics[scale=0.35]{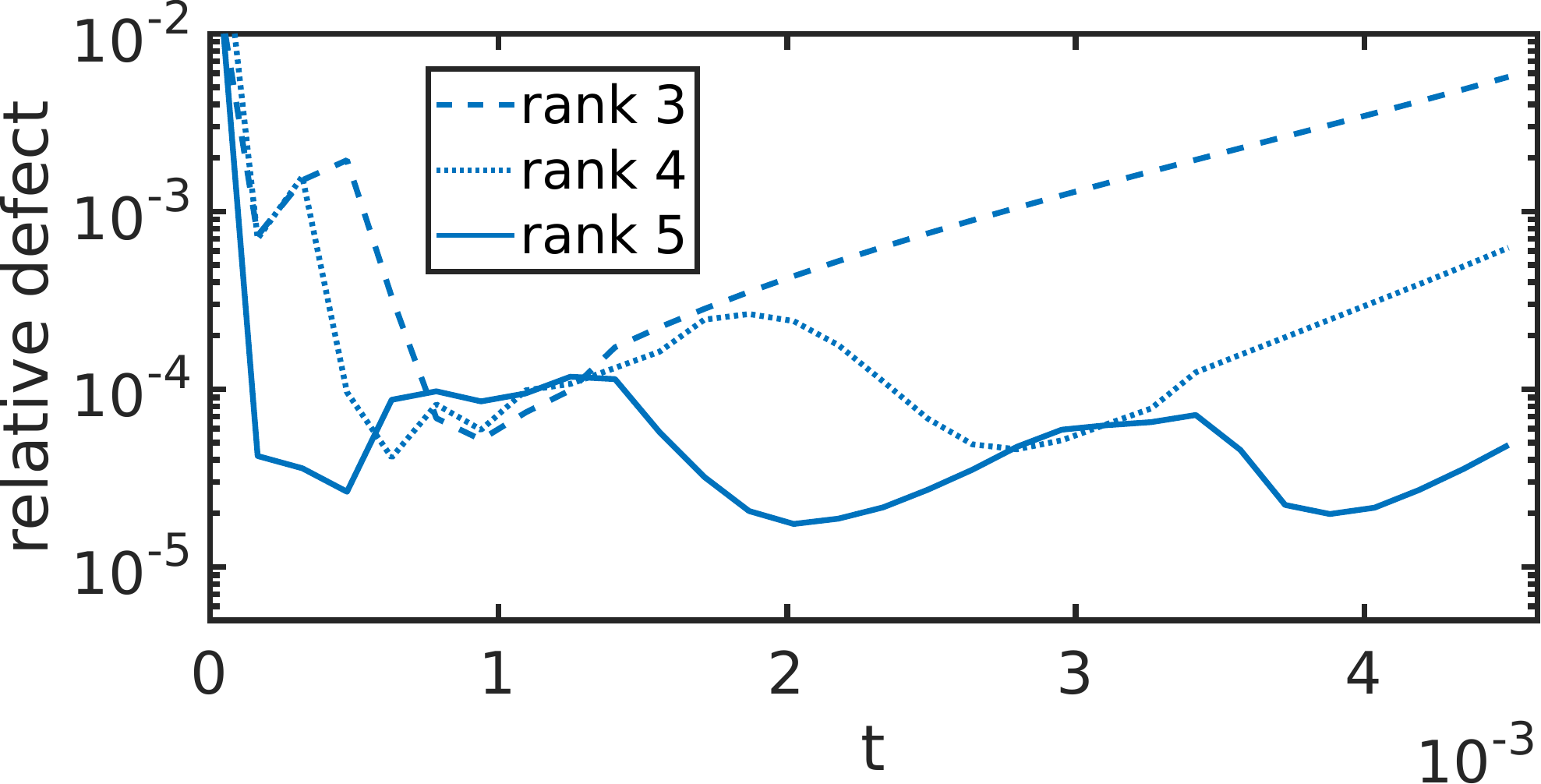}
\caption{Local relative defect of the dynamical approximation of \eqref{eqn:reacdiff} in the case $d=3$ and $I=400$ for different ranks $\mathbf{r} = (r,r,r)$.}\label{fig:reacdiff3d}
\end{figure}
Fig.~\ref{fig:timings} shows average timings for one iteration in Alg.~\ref{alg:alseuler} for varying mode size $I$ and ranks $r$ in the case $d=3$. The function evaluation involves compression (via HOOI) of Hadamard products 
with squared ranks, which is the most expensive part. For larger ranks the part of $\mathbf{B}_n$ involving the already computed right hand side scales with $r^3$. The rest, especially scalings with respect to 
the mode size, scales constant or at most linearly in the range of consideration.
\begin{figure}[hbtp]
\hspace{0.0 in} 
\centering 
\includegraphics[scale=0.35]{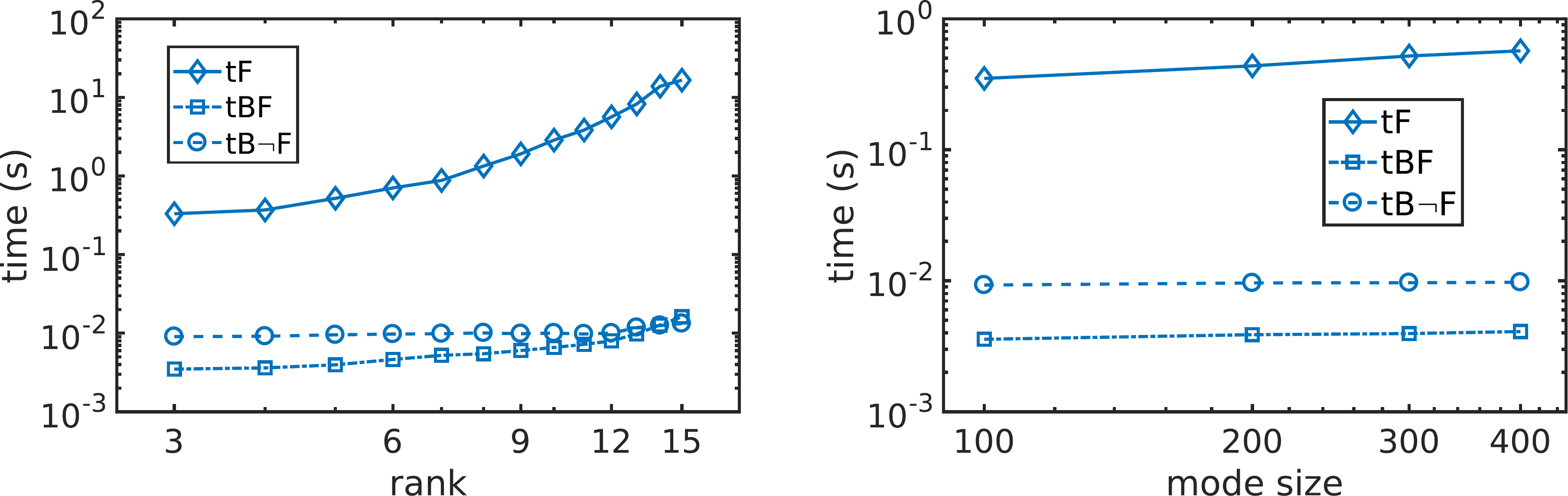}
\caption{Timings for one iteration in Alg.~\ref{alg:alseuler} for the discretization of problem \eqref{eqn:reacdiff} in the case $d=3$. Left: Varying rank 
$r$ for fixed mode size $I=300$. Right: Varying mode size $I$ for fixed rank $r=5$.  In each case timings is given for the right hand side evaluation (tF), the part of the $\mathbf{B}_n$ computation involving the already evaluated right hand side (tBF), 
as well as the part without right hand side and SVD (tB$\neg$F).}\label{fig:timings}
\end{figure}
}

\section{Conclusion}
An optimization-based approach for dynamical low-rank matrix and Tucker tensor approximation of parameter-dependent data tensors and solutions of tensor differential equations such as MCTDH of molecular quantum dynamics was presented. 
The approach is motivated by the effectiveness of alternating schemes for static tensor fitting problems, which reduce quadratic optimization subproblems to eigenproblems accessible via SVD. 
The presented method uses alternating least squares to establish an Euler scheme for the dynamical fitting problem on the tangent space without a gauge constraint. The quadratic subproblems are found to be 
equivalent to trace optimization with solutions being explicitly computable via SVD of small size. Remarkably, the method needs no pseudo inverse of matrix unfoldings of the core tensor 
(or inverse of density matrices in MCTDH), which enhances stability in the presence of small singular values. 
This makes the method robust and eliminates numerical stability problems coming from ''overestimated'' initial values{, which is confirmed by means of computational results.} 
Regularization of Tikhonov type is suggested to compensate for the lack of uniqueness of the tangent space representations, 
having the positive side effect to further stabilize in rank deficient cases. Higher order methods can be composed, where a second order example is given. 
In forthcoming work the approach can be extended to the hierarchical Tucker format and tensor trains (matrix product states) \cite{lubich2013dynamical}.

\section*{Acknowledgments}
Financial support by the Austrian Science Fund (FWF) via the SFB "ViCoM" (grant F41) {and SFB "Complexity in PDEs" (grant F65)} is acknowledged. 
The computations were achieved by using the Vienna Scientific Cluster (VSC).

\bibliography{bibref}

\end{document}